\documentclass[A4paper,11pt]{article}
\usepackage{latexsym}
\usepackage{mathrsfs}
\usepackage{amssymb}
\usepackage{amscd}
\usepackage[dvips]{graphicx}                  

\newtheorem{theorem}{Theorem}
\newtheorem{corollary}{Corollary}

\newenvironment{proof}[1][Proof]{\textbf{#1.} }{\ \rule{0.5em}{0.5em}}
\date{}
\long\def\symbolfootnote[#1]#2{\begingroup%
	\def\thefootnote{$\;$}\footnote[#1]{$^*$#2}\endgroup}
\begin{document}
	
	\title{Some remarks on polarized partition relations}
	\author{Joanna Jureczko}
\maketitle

\symbolfootnote[2]{Mathematics Subject Classification: Primary 03E02, 03E05, 03E10, 54E20.

\hspace{0.2cm}
Keywords: \textsl{partition relation, polarized  partition relation, strong sequences, product of generalized strong sequence, cardinal number, ordinal number, strongly inaccessible number}}

\begin{abstract}
	This paper deals with two notions: a polarized partition relations $\left( \begin{array}{c}
	\alpha \\
	\beta 
	\end{array} \right) \to \left( \begin{array}{cc}
	\gamma & \eta \\
	\delta & \lambda 
	\end{array} \right)$ and product of generalized strong sequences.  Strong sequences were introduced by Efimov in 1965 as a usefull tool for proving famous theorems in dyadic spaces, i.e. continuous images of Cantor cube. In this paper we introduce the notion of product of generalized strong sequences and give the pure combinatorial proof that $\left( \begin{array}{c}
	\alpha \\
	\beta 
	\end{array} \right) \to \left( \begin{array}{cc}
	\gamma & \eta \\
	\delta & \lambda 
	\end{array} \right)$
is a consequence of 	the existence of product of generalized strong sequences.
	\end{abstract}

\section{Introduction and historical background}

The notion of partition relations was introduced in \cite{ER} by Erd\"os and Rado as the ordinary partition relations which concerned partitions of finite subsets of a set of a given size and the polarized partition relations which concerned partitions of finite subsets of products of set of a given size, (where size means cardinality of a set or order type of an ordered set - we will specify it in the concrete situations). However, this topic has its origin in paper the Ramsey's paper \cite{R}. The main result of \cite{R} was generalized in 1942 by Erd\"os, \cite{E}.

Papers that deserve attention in this topic are undoubtedly \cite{ER2, HL, W}, however a great many new results were proved by researchsers in the recently time. This shows that the topic is extremely lively and still worth exploring.

In paper \cite{JJ1}, we proved a several theorems which we call Ramsey type  contain an alternative in the thesis: either we obtain a set of large cardinality with a certain property, or we obtain a set with small cardinality with an opposite property. This lead us to the following array notation
$$(\alpha) \to (\beta, \gamma)^n$$
which means that for given cardinals $\alpha, \beta, \gamma$ and for each set $A$ of cardinality $\alpha$ and a function $c \colon [A]^n \to 2$ there exists  a set $A_0 \subseteq A$ of cardinality $\beta$ such that $c|[A_0]^2 = \{0\}$ or there exists a set $A_1 \subseteq A$ of cardinality $\gamma$ such that $c|[A_1]^2 = \{1\}$.
In the literature there are known significant results of such theorems. We recall here some of them.
The result by Hajnal \cite{H} says that if $2^\kappa = \kappa^+$, then $$\kappa^+ \not \to (\kappa^+, \kappa+1)^2.$$ Todorcevic, in \cite{T}, proved that PFA (Proper Forcing Axiom) implies $$\omega_1 \to (\omega_1, \alpha)^2$$ for all $\alpha < \omega_1.$ Dushnik and Miller in \cite{DM} showed that for every infinite cardinal $\kappa$ $$\kappa \to (\kappa, \omega)^2.$$ 
Chang, in \cite{C}, proved that for all $m< \omega $
$$\omega^\omega \to (\omega^\omega, m)^2.$$

 Adopting the array notation to te statement in \cite{JJ1}, we obtain that for cardinals $\beta, \kappa, \tau $ such that $\omega \leqslant \beta \ll \tau, \kappa < \tau$ and $\tau, \beta$ regular 
 $$\tau \to (\tau, \beta),$$ where  $\beta \ll \tau$ denote $\tau$ is \textit{strongly $\beta$-inaccessible}, i.e $\beta < \tau$ and $\alpha^\lambda < \tau$ whenever $\alpha < \beta, \lambda < \tau$.  
 
 In the literature one can meet with "a kind of combination" of the above two partition relations, i.e.
  $$\left( \begin{array}{c}
  \alpha \\
  \beta 
  \end{array} \right) \to \left( \begin{array}{cc}
  \gamma & \eta \\
  \delta & \lambda 
  \end{array} \right)$$
  (for a given system of numbers $\alpha, \beta, \gamma, \delta, \eta, \lambda, \kappa$). For example Erd\"os, Rado and Hajnal, in \cite{ER2},  showed that
  $$\left( \begin{array}{c}
  \lambda \\
  \kappa 
  \end{array} \right) \to \left( \begin{array}{cc}
  \lambda & n \\
  \kappa & \kappa 
  \end{array} \right)$$
  holds whenever $n< \omega\leqslant \kappa < \lambda$ and $\lambda$
 regular. Jones in \cite{J} showed the direct proof of the following result: if $\kappa = \omega$ and $\lambda=2^{< \kappa}$ then
 $$\left( \begin{array}{c}
\lambda^+ \\
\lambda 
\end{array} \right) \to \left( \begin{array}{cc}
\lambda^+ & \alpha \\
\lambda & \kappa
\end{array} \right)$$
for all $\alpha < \omega_1.$ 

In \cite{W} we have more two interesting results. The first one is 
$$\left( \begin{array}{c}
\kappa \\
\kappa^+ 
\end{array} \right) \to \left( \begin{array}{cc}
\kappa & \eta \\
\kappa^+ & \kappa^+
\end{array} \right)$$
for any infinite $\kappa$, if $\eta< \kappa$.
While, the second is  
$$\left( \begin{array}{c}
\kappa \\
\kappa^+ 
\end{array} \right) \to \left( \begin{array}{cc}
\kappa & \kappa \\
\kappa & \eta
\end{array} \right)$$
for any $\kappa-$ singular and $\eta< \kappa$.
\\
  Some  newer interesting results in this topic one can find e.g. in \cite{KW}.
  
  The main result of this paper is to prove that the following theorem
   $$\left( \begin{array}{c}
  \alpha_1 \\
  \alpha_2\\
  ...\\
  \alpha_n 
  \end{array} \right) \to \left( \begin{array}{cc}
  \beta_1 & \gamma_1 \\
 \beta_2 & \gamma_2 \\
 ...\\
 \beta_n & \gamma_n \\
  \end{array} \right)$$
  is the consequence Theorem 1 in Section 3 which is also in \cite{JJ5}. 
  To make this work self-sufficient we cite the proof of theorem on product strong sequences from \cite{JJ5}.
  
  This paper is a continuation of \cite{JJ5} in which we show that the polarized partition relation
  $$\left( \begin{array}{c}
  \alpha_1 \\
  \alpha_2 \\
  ...\\
  \alpha_n 
  \end{array} \right) \to \left( \begin{array}{c}
  \beta_1 \\
  \beta_2 \\
  ...\\
  \beta_n 
  \end{array} \right)_\gamma$$
  is equivalent to the existence of strong sequences.
  
  The strong sequences method was introduced by Efimov in \cite{E3} as a useful tool for proving theorems in dyadic spaces, (i.e. continuous images of the Cantor cube). Among others, Efimov showed that strong sequences does not exist in the general Cantor discontinua.
  The  topic of strong sequences was considered by Turza\'nski in the 90s' of the last century. Turza\'nski reformulated the definition of strong sequences as follows. 
  
  Let $X$ be a set and let $B \subseteq P(X)$ be a family of non-empty subsets of $X$ closed with respect to finite intersections. Let $H_\alpha \subseteq B$ and $S_\alpha\subseteq B$ such that $S_\alpha$ is finite. A sequence $(S_\alpha, H_\alpha)_{\alpha < \kappa}$ is called a \textit{strong sequence} iff $S_\alpha \cup H_\alpha$ is centered and $S_\beta \cup H_\alpha$ is not centered, whenever $\beta > \alpha$.
  
  In \cite{MT1} the author proved the following theorem on strong sequences: if there exists a strong sequence   $(S_\alpha, H_\alpha)_{\alpha <(\kappa^\lambda)^+}$ such that $|H_\alpha| \leqslant \kappa$ for all $\alpha < \kappa^\lambda)^+$ then the family $B$ contains a subfamily of cadinality $\lambda^+$ consisting of pairwise disjoint sets. 
  
  Based on this result Turza\'nski estimated the weight of regular spaces.
  In \cite{MT2}  he gave a new proof of Esenin-Volpin Theorem of weight of dyadic spaces (in general form in class of thick spaces which possesses special subbases) and in \cite{MT3} he showed that the theorem on strong sequences is equivalent to Erd\"os-Rado Theorem. 
  
  The investigations on strong sequences have been continued, extended and improved by the author of this paper in \cite{JJ3, JJ4, JJ, JJ1, JJ2, JJ6, JJ5}. 
  
  In paper \cite{JJ} there is proved the generalization of theorem on strong sequences and it is shown that it is equivalent to the generalized Erd\"os-Rado Theorem. Further generalizations of these results are given in \cite{JJ1}. In \cite{JJ3, JJ6}   there is introduced the cardinal invariant associated with strong sequences and there are shown some inequalities between it and well known cardinal invariants. In \cite{JJ2} there is shown that the existence of so called K-Lusin sets is equivalent to the existence of strong sequences of the same cardinality.  
  The newest result concerning strong sequences are concentrated around product of generalized strong sequences and its connections with polarized partition relations, (\cite{JJ5}).
  However, we know a number of consequences of the existence of strong sequences the topic seems not to be exhausted. The main problem followed from Efimov result is still open: if strong sequences does not exist in general Cantor discontinua for which spaces does they exist?
  
  The paper is organized as follows. In Section 2 there are given basic definitions needed in further parts of the text. In section 3 there is proved the theorem on product of strong sequences. In Section 4 there is shown the equivalnece of Theorem 1 with polarized partition relation. 
  
  In this paper there are used standard notation and terminology. For the definitions and facts not cited her we refer the reader to \cite{E2, J2}.
  
  \section{Definitions}
  
  In the whole paper we use Greek letters to denote the cardinal or ordinal numbers, (which one we will mean at the particular parts will be follow from the context). 
  \\\\
  
  \textbf{2.1.} The \textit{polarized partition relation} 
  $$\left( \begin{array}{c}
  \alpha_1 \\
  \alpha_2 \\
  ... \\
  \alpha_n 
  \end{array} \right) \to \left( \begin{array}{cc}
  \beta_1 &\gamma_1\\
  \beta_2 &\gamma_2\\
  ... & ...\\
  \beta_n &\gamma_n
  \end{array} \right)$$
  means that the following  statement is true: for every set $A_k$ of cardinality $\alpha_k$, $(1 \leqslant k \leqslant n)$  and for every function $$c \colon A_1\times A_2 \times ... \times A_n \to 2$$ either there are $B_k \subseteq A_k, (1 \leqslant k \leqslant n)$ of cardinality $\beta_k$ such that $$c|(B_1\times B_2 \times ... \times B_n) = \{0\}$$ or there are $C_k \subseteq A_k, (1 \leqslant k \leqslant n)$ of cardinality $\gamma_k$ such that $$c|(C_1\times C_2 \times ... \times C_n) = \{1\}.$$ 
  \noindent
  \textbf{2.2.} Let $(X_k, r_k)$ be sets with two-place relations $r_k$, $(1\leqslant k \leqslant n)$.
  
  In the whole paper we restrict our considerations to finite products of sets, because we do not need more in this moment, but the results presented in further parts of this paper can be generalized for infinite products, (with extreme caution as is usual with infinite product operations). 
  
  Let $|X_k| = \kappa_k, (1\leqslant k \leqslant n)$. Let $$X = X_1\times X_2\times ... \times X_n$$ and $\kappa = \kappa_1\cdot \kappa_2 \cdot ... \cdot \kappa_n$.
  
  We say that $a = (a_1, a_2, ..., a_n), b = (b_1, b_2, ..., b_n) \in X$ \textit{have a bound} iff  there is $c= (c_1, c_2 ,..., c_n)$ such that $(a_k, c_k) \in r_k$ and $(b_k, c_k) \in r_k$ for every $1\leqslant k \leqslant n$.
  
  We say that $A \subseteq X$ is \textit{$\kappa$-directed} if every subset of $A$ of cardinality less than $\kappa$ has a bound.
  \\\\
  \noindent
  \textbf{2.3.} Let $(X_k, r_k)$ be sets with relations $r_k$ and $\alpha, \kappa_k, (1 \leqslant k \leqslant n)$, be cardinals. 
  A sequence $(H^k_\xi)_{\xi< \alpha}$ is called a \textit{$\kappa_k$-strong sequence} iff
  \begin{itemize}
  	\item [(1)] $H^k_\xi$ is $\kappa_k$-directed for all $\xi < \alpha$
  	\item [(2)] $H^k_\xi \cup H^k_\psi$ is not $\kappa_k$-directed whenever $\xi < \psi < \alpha$, i.e. there exists $S^k_\psi \in [H^k_\psi]^{<\kappa_k}$ such that for any $\xi < \psi$ the set $H^k_\xi \cup S^k_\psi$ is not $\kappa_k$-directed. (Such a set $S^k_\psi$ is called \textit{$(k, \xi, \psi)$-destroyer}.).
  \end{itemize}
  Denote $X = X_1\times X_2\times ... \times X_n$, $\kappa = \kappa_1 \cdot \kappa_2 \cdot ... \cdot \kappa_n$ and  $H_\xi= H^1_\xi \times H^2_\xi \times ... \times H^n_\xi$.
  A sequence $(H_\xi)_{\xi< \alpha}$ is called a \textit{product $\kappa$-strong sequence} iff 
  \begin{itemize}
  	\item [(3)] $H_\xi$ is $\kappa$-directed for all $\xi < \alpha$
  	\item [(4)] $H_\xi \cup H_\psi$ is not $\kappa$-directed whenever $\xi < \psi < \alpha$,
  \end{itemize}

\noindent
\textbf{2.4.}  An uncountable cardinal $\tau$ is \textit{weakly inaccessible} iff it is a limit cardinal and is regular.
A cardinal $\tau$ is \textit{(strongly) inaccessible} if $\tau$ is an infinite regular cardinal and such that $2^\beta < \tau$, whenever $\beta < \tau$. A cardinal $\tau$ is \textit{strongly $\beta$-inaccessible} iff $\beta < \tau$ and $\alpha^\lambda < \tau$, whenever $\alpha < \tau, \lambda < \beta$. By $\beta \ll \tau$ we denote that $\tau$ is strongly $\beta$-inaccessible.

Every inaccessible cardinal is weakly inaccessible. If the GCH holds then every weakly inaccessible cardinal $\tau$ is inaccessible.
Inaccessible cardinals cannot be obtained from smaller cardinals by the usual set-theoretical operations. This is one of the themes of set theory of large cardinals. 
The existence of inaccessible cardinals is not provable in ZFC. Moreover, it cannot be shown that the existence of inaccessible cardinals is consistent with ZFC, (see \cite[Theorem 12.12]{J2}).
The least inaccessible cardinal is not measurable, (see \cite{WH}).
Inaccessible cardinals were introduced by Sierpi\'nski and Tarski in 1930.

  \section{Theorem on product $\kappa$-strong sequences}

  The special case of the following theorems (for $k=1$) was proved in \cite{JJ1}.
  \\
  
\begin{theorem}  Let $n < \omega$. For  $1\leqslant k \leqslant n$ let $\beta_k, \kappa_k, \mu_k, \eta_k$ be cardinals such that $\omega\leqslant \beta_k \ll\eta_k, \mu_k < \beta_k, \kappa_k \leqslant 2^{\mu_k}, \eta = \eta_1\cdot \eta_2 \cdot ...\cdot \eta_n, \kappa = \kappa_1\cdot \kappa_2 \cdot ...\cdot \kappa_n$  and $\beta_k, \eta_k$ be regulars. Let $X= X_1\times X_2\times ... \times X_n$ be a set and $|X_k| = \eta_k$. If there exists a product $\kappa$-strong sequence $$\{H_\alpha \subseteq X \colon \alpha < \eta\}$$ with $|H^k_\alpha| \leqslant 2^{\mu_k}$ for all $\alpha < \eta$, then there exists a product $\kappa$-strong sequence $$\{T_\alpha \colon \alpha < \beta\}$$ with $|T^k_\alpha|<\kappa_k$ for all $\alpha < \beta$.
\end{theorem}

\begin{proof}
	Fix $n < \omega$.
	Let $\{H_\alpha \colon \alpha < \eta\}$  be such that $H_\alpha = H^1_\alpha \times H^2_\alpha \times ... \times H^n_\alpha$ is a $\kappa$-strong sequence. Fix $\alpha$ and name it $\alpha_0$, (without the loss of generality one can assume that $\alpha_0 = 0$).
	For each $k, (1 \leqslant k \leqslant n)$, consider a function 
	$$f^k_{\alpha_0}\colon \eta_k\setminus \{\alpha_0\} \to [H^k_{\alpha_0}]^{< \kappa_k}$$
	such that $f^k_{\alpha_0}(\xi)=S^k_{\alpha_0}$ for some $S^k_{\alpha_0} \in [H^k_{\alpha_0}]^{< \kappa_k}$. Since $|H^k_{\alpha_0}| \leqslant 2^{\mu_k}$, hence $|[H^k_{\alpha_0}]^{< \kappa_k}|\leqslant 2^{\mu_k}< \eta_k$. It means that $f^k_{\alpha_0}$ determines a partition of $\eta_k\setminus \{\alpha_0\}$ into at most $2^{\mu_k}< \eta_k$ elements. The cardinal $\eta_k$ is regular, hence there exists $\overline{S}^k_{\alpha_0}\in [H^k_{\alpha_0}]^{< \kappa_k}$ such that $|(f^k_{\alpha_0})^{-1}(\overline{S}^k_{\alpha_0})|= \eta_k$.
	Let $$\mathcal{S}^k_{\alpha_0} = \{S^k_{\alpha_0} \in [H^k_{\alpha_0}]^{< \kappa_k} \colon |(f^k_{\alpha_0})^{-1}(S^k_{\alpha_0})|= \eta_k\}$$
	and let
	$$\mathcal{A}^k_{\alpha_0} = \{A^k_{\alpha_0}\subseteq \eta_k\setminus \{\alpha_0\} \colon \exists(S^k_{\alpha_0}\in [H^k_{\alpha_0}]^{< \kappa_k})\  (f^k_{\alpha_0})^{-1}(S^k_{\alpha_0})= A^k_{\alpha_0} \textrm{ and } |A^k_{\alpha_0}| = \eta_k\}$$
	be a family of pairwise disjoint sets. 
	
	Before the continuation of the proof we make the following observation. By definition of product $\kappa$-strong sequence for each $\alpha> \alpha_0$ there exists a $(k, \alpha_0, \alpha)$-destroyer. Since we consider only $k$ such that $(1\leqslant k \leqslant n)$ some of them must occur at least $\beta_k$-times, $(\omega \leqslant \beta_k \ll \eta _k)$.
	Now, we are ready to continue the proof.	
	
	For  every relevant $k$ we will construct inductively 
	\begin{itemize}
		\item [a)] an increasing subsequence $\{\alpha_\gamma \colon \gamma < \beta_k\}$ of elements from  $\eta_k,$
		\item [b)] families $\mathcal{A}^k_{\alpha_\gamma} = \{A^k_{\alpha_\gamma}\subseteq \eta_k \colon \exists(S\in [H^k_{\alpha_\gamma}]^{< \kappa_k})\  (f^k_{\alpha_\gamma})^{-1}(S)= A^k_{\alpha_\gamma} \textrm{ and } |A^k_{\alpha_\gamma}| = \eta_k\} $,
	\end{itemize}
	where 	$$f^k_{\alpha_\gamma}\colon A \setminus \{\alpha_\gamma\}\to [H^k_{\alpha_\gamma}]^{< \kappa_k}$$
	such that $A \in \mathcal{A}^k_{\alpha_{\gamma-1}}$ and $f^k_{\alpha_\gamma}(\gamma) = S$ for some $(k, \alpha_\sigma, \alpha_\gamma)$-destroyer, whenever $\alpha_\sigma < \alpha_\gamma$.
	
	Assume that we have constructed increasing subsequence $\{\alpha_\gamma \colon \gamma < \beta_k\}$ of $\eta_k$ and families $\mathcal {A}^k_{\alpha_\gamma}$ as was done above. 
	
	Next, choose $\alpha > \alpha_\gamma, (\alpha < \eta_k)$ such that there exists a $(k, \alpha_\gamma, \alpha)$-destroyer $S^k_{\alpha} \in [H^k_{\alpha}]^{< \kappa_k}$ and denote this $\alpha$ by $\alpha_\delta$, where $\delta= \gamma+1$.
	For each $A^k_{\alpha_\gamma} \in \mathcal{A}^k_{\alpha_\gamma}$ define a function 
	$$f^k_{\alpha_\delta}\colon A^k_{\alpha_\gamma} \setminus \{\alpha_\delta\}\to [H^k_{\alpha_\delta}]^{< \kappa_k}$$
	such that $f^k_{\alpha_\delta}(\xi) = S^k_{\alpha_\delta}$ for some $(k, \alpha_\gamma, \alpha_\delta)$-destroyer $S^k_{\alpha_\delta} \in [H^k_{\alpha_\delta}]^{\kappa_k}$ (and any $\alpha_\gamma < \alpha_\delta$).
	Since $|H^k_{\alpha_\delta}| \leqslant 2^{\mu_k}< \eta_k$ and $|A^k_{\alpha_\gamma}| = \eta_k$ the function $f^k_{\alpha_\delta}$ determines a partition of $A^k_{\alpha_\gamma} \setminus \{\alpha_\delta\}$ into at most $2^{\mu_k}$ elements. Hence there exists a $(k, \alpha_\gamma, \alpha_\delta)$-destroyer $\overline{S}^k_{\alpha_\delta} \in [H^k_{\alpha_\delta}]^{< \kappa_k}$ such that $|(f^k_{\alpha_\delta})^{-1}(\overline{S}^k_{\alpha_\delta})|= \eta_k$. Let 
	$$\mathcal{S}^k_{\alpha_\delta} = \{S^k_{\alpha_\delta} \in [H^k_{\alpha_\delta}]^{< \kappa_k} \colon |(f^k_{\alpha_\delta})^{-1}(S^k_{\alpha_\delta})|= \eta_k\}$$
	and let
	$$\mathcal{A}^k_{\alpha_\delta} = \{A^k_{\alpha_\delta}\subseteq  A^k_{\alpha_\gamma} \setminus \{\alpha_\delta\}\colon \exists_{S^k_{\alpha_\delta}\in [H^k_{\alpha_\delta}]^{< \kappa}} |(f^k_{\alpha_\delta})^{-1}(S^k_{\alpha_\delta})| =\eta_k\}$$
	be a family of pairwise disjoint sets.
	
	If $\delta$ is limit, we consider 
	$$f^k_{\alpha_\delta}\colon \bigcap_{\rho< \gamma }A^k_{\alpha_\rho} \setminus \{\alpha_\delta\}\to [H^k_{\alpha_\delta}]^{< \kappa_k}$$
	for $\bigcap_{\rho< \gamma }A^k_{\alpha_\rho} \not = \emptyset$.
	The induction step is complete.
	
	Now define a sequence 
	$$T_{\alpha_\gamma} = T^1_{\alpha_\gamma} \times T^2_{\alpha_\gamma} \times ... \times T^n_{\alpha_\gamma}$$ for all $\alpha_\gamma< \beta = \beta_1\cdot \beta_2\cdot ... \cdot \beta_n$ in the following way:
	\[ T^k_{\alpha_\gamma}=
	\left\{\begin{array}{rl}
	S^k_{\alpha_\gamma} \in  [H^k_{\alpha_\gamma}]^{< \kappa_k} & \mbox{if }\exists k\  S^k_{\alpha_\gamma} \in \mathcal{S}^k_{\alpha_\gamma}\\
	P^k_{\alpha_\gamma} \in  [H^k_{\alpha_\gamma}]^{< \kappa_k} & \mbox{otherwise},
	\end{array}\right. \]
	where $P^k_{\alpha_\gamma}$ is chosen arbitrary. Thus, we have defined at least one product $\kappa$-strong sequence $\{T_{\alpha_\gamma} \colon \alpha_ \gamma < \beta\}$ of the required property.
	
	Suppose now, that at least one of product $\kappa$-strong sequence $$\{T_{\alpha_\gamma} \colon \alpha_ \gamma < \beta\}$$ has length $\zeta> \beta$, i.e. there exists $k, (1\leqslant k \leqslant n)$ which occurs $\zeta > \beta$-times, i.e. there is a sequence $$\{S^k_{\alpha_\gamma} \colon \alpha_ \gamma < \zeta\}$$
	such that $S^k_{\alpha_\gamma} \in \mathcal{S}^k_{\alpha_\gamma}$.
	By our construction each $S^k_{\alpha_\gamma}$ determines a set $A \in \mathcal{A}^k_{\alpha_\gamma}$ such that $|A| = \eta_k$.
	Let $$\nu_k = \sup\{|A|\colon \alpha_\gamma, A \in S^k_{\alpha_\gamma}\}.$$
	Then, there would exist $\nu_k^\zeta > \eta_k$ pairwise disjoint sets $A \in \mathcal{A}^k_{\alpha_\gamma},$ where $ |A| = \eta_k$. A contradiction.
\end{proof}

\begin{theorem}  Let $n < \omega$. For  $1\leqslant k \leqslant n$ let $\beta_k, \kappa_k, \mu_k, \eta_k$ be cardinals such that $\omega\leqslant \beta_k \ll\eta_k, \mu_k < \beta_k, \kappa_k \leqslant 2^{\mu_k}, \eta = \eta_1\cdot \eta_2 \cdot ...\cdot \eta_n, \kappa = \kappa_1\cdot \kappa_2 \cdot ...\cdot \kappa_n$  and $\beta_k, \eta_k$ be regulars. Then
	either $X$ contains a $\kappa$-directed subset of cardinality $\eta$ or $X$ contains a subset of cardinality $\beta$ consisting of elements which are pairwise disjoint.
	\end{theorem}

\begin{proof}
	Without the loss of generality we can assume that $X_k\subseteq \eta_k, (1\leqslant k \leqslant n)$.
	Suppose that each $\kappa_k$-directed subset of $X_k$  has cardinality less than $\eta_k$. We will use Theorem 1 to prove the second alternative of the statement. 
	
	In order to do this we will construct inductively a product $\kappa$-strong sequence $$\{H_\alpha \colon \alpha < \eta\}.$$
	
	For each $k, (1\leqslant k \leqslant n)$ choose a $\kappa_k-$directed set $H^k_0 \subseteq X_k$. It is the first step of our induction.
	
	Assume that for $\alpha< \eta$ the product $\kappa$-strong sequence $$\{H_\gamma \subseteq X\setminus \bigcup_{\delta< \gamma} H_\delta \colon \gamma < \alpha\}$$ has been defined.
	Since $|H^k_\gamma|< \eta_k$ for each $k$ then $|\bigcup_{\gamma< \alpha} H^k_\gamma|< \eta_k$ and $\eta_k$ are regular we have $|X_k \setminus(\bigcup_{\gamma< \alpha} H^k_\gamma)|= \eta_k$.
	Hence we can construct $H^k_\alpha.$
	
	If $\alpha = \gamma+1$ then take a maximal $\kappa_k$-directed set 
	$$H^k_\alpha \subseteq X_k \setminus \bigcup_{\gamma< \alpha} H^k_\gamma$$ such that $H^k_\alpha\cup H^k_\gamma$ is not $\kappa_k$-directed for some $k$.
	
	If $\alpha$ is limit, then take $$x_k = \min(X_k \setminus \bigcup_{\gamma < \alpha} H^k_\gamma)$$ and put $$H^k_\alpha = \bigcup_{\gamma < \alpha} \cup \{x_k\}.$$ Obviously $H^k_\alpha \cup H^k_\gamma$ is not $\kappa_k$-directed for some $k$.
	
	Take $$H_\alpha = H^1_\alpha \times H^2_\alpha \times ... \times H^n_\alpha$$  as the next element of product $\kappa$-strong sequence.
	By Theorem 1, there exists a product $\kappa$-strong sequence
	$$\{T_\alpha = T^1_\alpha \times T^2_\alpha \times ... \times T^n_\alpha \colon \alpha < \beta\}$$ such that $T^k_\alpha \subseteq H^k_\alpha $ and $ |T^k_\alpha|< \kappa_k$ for all $\alpha < \beta$.
	If $T^k_\alpha, (\alpha < \beta_k)$ are not pairwise disjoint then take
	$$(T^k_0)'=T^k_0,$$ $$(T^k_\alpha)' = T^k_\alpha \setminus \bigcup_{\gamma < \alpha} T^k_\gamma \textrm{ for } \alpha-\textrm{successor}$$ and $$(T^k_\alpha)' = \sup (\bigcup_{\gamma < \alpha}(T^k_\gamma)') \textrm{ for } \alpha-\textrm{limit}.$$	
\end{proof}
\\

As Corollary of Theorem 2 we obtain  
  
\begin{corollary}  Let $n < \omega$. For  $1\leqslant k \leqslant n$ let $\beta_k, \kappa_k, \mu_k, \eta_k$ be cardinals such that $\omega\leqslant \beta_k \ll\eta_k, \mu_k < \beta_k, \kappa_k \leqslant 2^{\mu_k}, \eta = \eta_1\cdot \eta_2 \cdot ...\cdot \eta_n, \kappa = \kappa_1\cdot \kappa_2 \cdot ...\cdot \kappa_n$  and $\beta_k, \eta_k$ be regulars. Then
	$$\left( \begin{array}{c}
	\eta_1 \\
	\eta_2 \\
	...\\
	\eta_n 
	\end{array} \right) \to \left( \begin{array}{cc}
	\eta_1 &\beta_1\\
	\eta_2 &\beta_2\\
	... & ...\\
	\eta_n &\beta_n
	\end{array} \right).$$
\end{corollary}

\begin{proof}
	Take 
for every set $A_k$ of cardinality $\eta_k$, $(1 \leqslant k \leqslant n)$  and take a colouring function $$c \colon A_1\times A_2 \times ... \times A_n \to 2.$$ If there are $B_k \subseteq A_k$ of cardinality $\kappa_k, (1 \leqslant k \leqslant n)$ such that the set $B_1\times B_2\times ... \times B_n$  is $\kappa$-directed then put
$$c|(B_1\times B_2 \times ... \times B_n) = \{0\}.$$
If there are $C_k \subseteq A_k$ of cardinality $\beta_k, (1 \leqslant k \leqslant n)$ such that the set $C_1\times C_2\times ... \times C_n$  consisting of elements which are pairwise disjoint then put
$$c|(C_1\times C_2 \times ... \times C_n) = \{1\}.$$
\end{proof}
\\
\\
\textbf{Declaration} The author have no conflicts of interest to declare.
\\\\
\textbf{Acknowledgments} The author would like to thank the reviewer for a detailed and insightful reading of the text and for all the comments that undoubtedly helped to improve the text and avoid inaccuracies and omissions.

\noindent
{\sc Joanna Jureczko}
\\
Wroc\l{}aw University of Science and Technology, Wroc\l{}aw, Poland
\\
{\sl e-mail: joanna.jureczko@pwr.edu.pl}

\end{document}